\documentclass{article}
\usepackage{amsmath,amssymb,bm}
\usepackage{amssymb,bm}
\usepackage{graphicx}
\usepackage{color}
\usepackage{amscd}

\usepackage{amsthm}

\newtheorem{theorem}{Theorem}

\newtheorem{prop}[theorem]{Proposition}

\newtheorem{remark}{Remark}

\newenvironment{pr}{\proof[\proofname]}{\endproof}

\begin{document}
 
\title
{Tangential morphisms via log arithmetic geometry}

\author
{Yuichiro Hoshi, Makoto Matsumoto, Chikara Nakayama}
\maketitle

\begin{abstract}
\noindent We give a reformulation of tangential morphisms (which is a generalization of Deligne's tangential base point) via log geometry. 
\end{abstract}

\section*{Introduction}
\label{s:intro}
  In this short article, we give a reformulation of tangential morphisms 
(which is introduced by the second author 
as a higher dimensional generalization of 
Deligne's tangential base point) via log arithmetic geometry. 
  We explain that in terms of log geometry, the definition of tangential morphisms become simple and their basic properties are easily deduced. 
  All contents must be well-known among experts. 
  But since it seems there are no literatures and it is hindering exchange among researchers, we would like to write them explicitly. 

  The log fundamental group is first introduced in \cite{FK}. 
  For the basic properties used in this article, see also \cite{I}, \cite{S}, \cite{V}, \cite{Ho}, and \cite{Kato:FI2}. 
  For a pro-fs log scheme, its log fundamental group is defined by the limit of Galois categories.

\section{Tangential morphisms}
\label{s:def}
  Let $\tau\colon Z \to V$ be a tangential morphism. 
  Let the notation be as in \cite{mat} 2.1.
  Endow $V^*$ with the fs log structure defined by $D_1,\ldots, D_l$, and endow $U$ with the pullback log structure. 
  Take a chart $\bold N^l \to {\cal O}_U$ on $U$ defined by $t_1,\ldots,t_l$.
  Let $\tilde Z$ be a pro-fs log scheme 
$Z\otimes_{\bold Z[\bold N^l]}\bold Z[({\bold Q}'_{\ge0})^l],$ where 
$Z$ is endowed with the pullback log structure from $V^*$ and ${\bold Q}'_{\ge0}$ means $\bigcup_{p\nmid n} \frac1n\bold N$. 
  Then we have a diagram 
$$Z \to \tilde Z^{\mathrm {cl}}  \gets \tilde Z \to V^*$$ of log schemes, where $Z$ in this time is regarded as a log scheme with the trivial log structure and $\tilde Z^{\mathrm {cl}}$ denotes the underlying pro-scheme of $\tilde Z$ with the trivial log structure. 

  We explain that the first two morphisms in this diagram 
$Z \to \tilde Z^{\mathrm {cl}}  \gets \tilde Z$
induce an isomorphism between the usual $\pi_1$ of $Z$ and $\pi_1^{\log}$ of $\tilde Z$. 
  First, the morphism 
$Z \to \tilde Z^{\mathrm {cl}}$ 
is regarded as a closed nil-immersion of the usual schemes 
and hence it induces the isomorphism of the usual $\pi_1$. 
  Next, any kummer \'etale (=k\'et) covering of $\tilde Z$ is strict \'etale 
because in general a k\'et covering is k\'et locally strict (cf.\ \cite{Kato:FI2} 
Proposition 2.7).
  Hence the morphism $\tilde Z \to \tilde Z^{\mathrm {cl}}$ induces an isomorphism of 
log $\pi_1$. 

  Since $\pi_1^{\log}$ of $V^*$ is canonically isomorphic to the tame $\pi_1$ of $(V,\bigcup_i D_i)$ 
(\cite{Ho} Proposition B.7), the morphisms induce a homomorphism of groupoids from $\pi_1$ of $Z$ to the tame $\pi_1$ of $V$. 

\begin{prop}
  The last homomorphism coincides with the one in {\rm \cite{mat}} Proposition $2.1$. 
\end{prop}

\begin{pr}
  We use the same notation in \cite{mat} 2.1.
  Let $W \to V$ be a finite \'etale cover of $V$ which is tamely ramified at 
$\bigcup_i D_i$.
  Then the $W^* \to V^*$ in \cite{mat} 2.1 is regarded as the corresponding finite k\'et covering, where $W^*$ is endowed with the natural log structure. 
  It is enough to show that the pullback of $W^*$ over $\tilde Z$ coincides with the pullback of $\tau^*_{fet}(W)$ defined in \cite{mat} 2.1. 
  But the morphism $\tilde Z \to V^*$ factors as $\tilde Z \to Z\{t\}\to V^*$ with the first arrow being strict, 
where $Z\{t\}$ is regarded as a pro-fs log scheme endowed with the log structure defined by 
$(\bold Q'_{\ge0})^l \to A_{\{t\}}; (a_i)_i \mapsto (t_i^{a_i})_i$.
  The pullback of $W^*$ over $Z\{t\}$ (in the category of saturated log schemes) is nothing but $N(W\vert_{Z\{t\}})$ endowed with the pullback log structure from $Z\{t\}$. 
  Hence its further pullback over $\tilde Z$ coincides with the pullback of $\tau^*_{fet}(W)$ as desired.
\end{pr}

\section{Homotopy}
In the notation in \cite{mat} 2.2, 
an infinitesimal homotopy from $\tau$ to $\tau'$ induces an isomorphism over $Z$ between 
$\tilde Z$ by $(t_i)$ and $\tilde Z$ by $(s_i)$.
  Hence, via the construction in Section \ref{s:def}, 
we have the following. 

\begin{prop}
  An infinitesimal homotopy yields a homotopy   
from $\tau^*_{fet}$ to $\tau'{}^*_{fet}$. 
\end{prop}

It is clear that this homotopy 
coincides with that in {\rm \cite{mat}} Proposition $2.2$.

\begin{remark}
 If $t_i = s_iu_i$ with some units $u_i \equiv 1$ modulo the defining ideal of $Z$ in $V^*$, $(t_i)$ and $(s_i)$ define the same chart on $Z$ (endowed with the pullback log structure from $V^*$).  
 Hence, $\tilde Z$ is also the same. 
 Thus the tangential morphisms essentially depend only on the choice of the basis of normal cotangent bundle of $Z$ in $V^*$. 
 The second author heard from Deligne that this is the origin of the terminology of the tangential base point. 
\end{remark}

\section{Composition}
Let the notation be as in \cite{mat} 2.3.

We have a commutative diagram 
$$\begin{CD}
\tilde Z'_{(t,s)} @>>> \tilde Z^*_{(t)} @>>> V^{**} \\
@VVV @VVV @. \\
\tilde Z'_{(s)} @>>> Z^* \\
@VVV @. @. \\
Z', 
\end{CD}
$$
where $V^{**}$ is endowed with the log by $(t_i)_i$ and $(s_j)_j$, $Z^*$ is endowed with the log only by $(s_j)_j$, 
$\tilde Z^*_{(t)}$ is $Z^*\otimes_{\bold Z[\bold N^l]}\bold Z[({\bold Q}'_{\ge0})^l]$ (where $Z^*$ is endowed with the pullback log structure from $V^{**}$), 
and 
$\tilde Z'_{(t,s)}$ (resp.\ $\tilde Z'_{(s)}$) is the pro-fs log scheme obtained by applying the construction in Section \ref{s:def} to $Z' \to V^{**}$ (resp.\ 
$Z' \to Z^*$) and the chart defined by both $(t_i)_i$ and $(s_j)_j$ (resp.\ only by $(s_j)_j$). 
Since $\tau'^*_{fet}$, $\tau^*_{fet}$, and $(\tau\tau')^*_{fet}$ in \cite{mat} 2.3 are defined by using $\tilde Z'_{(s)}$, $\tilde Z^*_{(t)}$, and $\tilde Z'_{(t,s)}$, respectively (by Proposition 1 for 
$\tau'^*_{fet}$ and $(\tau\tau')^*_{fet}$, and by a similar argument for $\tau^*_{fet}$), we have the following.

\begin{prop}
There is a canonical homotopy between $\tau'^*_{fet} \circ \tau^*_{fet}$ and $(\tau\tau')^*_{fet}$. 
\end{prop}

It is straightforward to see that this homotopy 
coincides with that in {\rm \cite{mat}} Proposition $2.3$.

\section{GAGA}
Let the notation be as in the case of a geometric tangential point 
$\theta \colon \mathrm{Spec}\,\bar k \to \bar V$
in \cite{mat} 2.4.2. 
Consider the associated analytic tangential point 
$Z:=\mathrm{Spec}\,{\bold C} \to V^{an}$ and 
$Z \to \tilde Z^{\mathrm {cl}}  \gets \tilde Z \to (V^*)^{an}$, where $\tilde Z$ is a pro-fs log analytic space constructed similarly as in Section \ref{s:def}.
Taking $(\cdot)^{\log}$, we have 
$Z \to \tilde Z^{\mathrm {cl}}  \gets (\tilde Z)^{\log} \to ((V^{*})^{an})^{\log}$.
  Let $\bold R^l$ be the universal cover of $Z^{\log}$ 
($Z$ here is endowed with the pullback log structure from $\bar V^*$). 
  Since the morphism $(\tilde{Z})^{\log} \to Z^{\log}$ induced by the natural morphism $\tilde{Z} \to Z$ is a universal profinite cover of $Z^{\log}$, there 
is a canonical map $\bold R^l \to (\tilde Z)^{\log}$. 
  Then the fiber functor $Q$ on $V^{an}$ defined by the pullback to $\bold R^l$ and taking the global section $\Gamma(\bold R^l, \cdot)$ is compatible with the tangential point $F_Q$ defined by $\theta$.
  Hence we have a homomorphism $\pi_1^{top}(((V^{*})^{an})^{\log};Q,R) \to \pi_1(V; F_Q, F_R)$, where $R$ is another point (tangential or usual). 

\begin{prop}
  The last homomorphism becomes an isomorphism after pro-finitely completing the left set. 
\end{prop}

\begin{pr}
  Replacing a tangential base point with a usual point, 
we may assume that both $Q$ and $R$ are usual points. 
  Then replacing $V^*$ with $V$, we reduce to the nonlog case in \cite{SGA1} (mentioned in \cite{mat} 2.4.2). 
\end{pr}

\vspace{0.5cm}
\noindent
{\bf Acknowledgments.}\ 
The authors thank J.\ C.\ for leading them to join this work.
Y.\ H.\ was supported by JSPS KAKENHI Grant Number 24K06668
and by the Research Institute for Mathematical Sciences, an
International Joint Usage/Research Center located in Kyoto University,
as well as the International Center for Research in Next Generation
Geometry.
C.\ N.\ was partly supported by the JSPS KAKENHI Grant Number 21K03199 and 26K06723.

\noindent Yuichiro Hoshi \\
Research Institute for Mathematical Sciences \\
Kyoto University, Kyoto 606-8502 \\
Japan

\noindent yuichiro@kurims.kyoto-u.ac.jp

\par\bigskip\par
\noindent Makoto Matsumoto

\noindent AMAGAERU Institute of Free Mathematics \\
2-37-6 Narita-Higashi, Suginami-ku, Tokyo 166-0015 \\ 
Japan

\noindent matmotmak@gmail.com

\par\bigskip\par
\noindent Chikara Nakayama

\noindent Department of Economics \\ Hitotsubashi University \\
2-1 Naka, Kunitachi, Tokyo 186-8601 \\ Japan

\noindent c.nakayama@r.hit-u.ac.jp

\end{document}